\documentclass[]{interact}

\usepackage{epstopdf}
\usepackage[caption=false]{subfig}

\usepackage[numbers,sort&compress]{natbib}
\bibpunct[, ]{[}{]}{,}{n}{,}{,}

\theoremstyle{plain}
\newtheorem{thm}{Theorem}[section]
\newtheorem{lemma}[thm]{Lemma}
\newtheorem{cor}[thm]{Corollary}

\theoremstyle{definition}

\theoremstyle{remark}

\usepackage{color}
\usepackage{algorithmic}
\usepackage{algorithm}
\usepackage{pdflscape}
\usepackage{multirow}

\DeclareMathOperator*{\argmin}{arg\,min}

\definecolor{pinegreen}{rgb}{0.0, 0.47, 0.44}

\def\pf{\medskip\noindent {\bf Proof.}~~}

\begin{document}


\title{Integer Linear Programming Formulations for Double Roman Domination Problem}

\author{
\name{Qingqiong Cai\textsuperscript{a}, Neng Fan\textsuperscript{b}\thanks{N. Fan: The corresponding author. Email: nfan@email.arizona.edu; Phone: +1(520) 621 6557},  Yongtang Shi\textsuperscript{c}, Shunyu Yao\textsuperscript{c}}
\affil{\textsuperscript{a}College of Computer Science,  Nankai University, Tianjin, China;\\ \textsuperscript{b}Department of Systems \& Industrial Engineering, University of Arizona, Tucson, AZ, USA;\\ \textsuperscript{c}Center for Combinatorics and LPMC, Nankai University, Tianjin, China}
}

\maketitle

\begin{abstract}
For a graph $G= (V,E)$, a double Roman dominating function (DRDF) is a function $f : V \to \{0,1,2,3\}$ having the property that if $f (v) = 0$, then vertex $v$ must have at least two neighbors assigned $2$ under $f$ or {at least} one neighbor $u$ with $f (u) = 3$, and if $f (v) = 1$, then vertex $v$ must have at least one neighbor $u$ with $f (u) \ge 2$. In this paper, we consider the double Roman domination problem, which is an optimization problem of finding the DRDF $f$ such that $\sum_{v\in V} f (v)$ is minimum. We propose {five integer linear programming (ILP) formulations and one mixed integer linear programming formulation with polynomial number of constraints for this problem. Some additional valid inequalities and bounds are also proposed for some of these formulations.} Further, we prove that {the first four models indeed solve the double Roman domination problem, and the last two models} are equivalent to the others regardless of the variable relaxation or usage of a smaller number of constraints and variables. Additionally, we use one ILP formulation to give an $H(2(\Delta+1))$-approximation algorithm. All proposed formulations and approximation algorithm are evaluated on randomly generated graphs to compare the performance.
\end{abstract}

\section{Introduction}
Throughout this paper, all graphs considered are finite, simple and undirected. Let $G=(V,E)$ be a graph with vertex set $V$ and edge set $E$. For every vertex $v \in V$, the open neighborhood of $v$ is the set $N(v)=\{u \in V \ | \ uv \in E\}$, and the closed neighborhood of $v$ is the set $N[v] = N(v) \cup \{v\}$. The degree of a vertex $v \in V$ is defined as $d_{G}(v) = |N(v)|$. The minimum and maximum degree of a graph $G$ are denoted by $\delta$ and $\Delta$, respectively. The distance of two vertices is the length of a shortest path between them in a graph. The diameter $d$ of a graph is the largest distance between vertices. The complement $\overline{G}$ of $G$ is the graph whose vertex set is $V$ and whose edges are the pairs of nonadjacent vertices of $G$. {The girth of $G$, denoted by $g$, is the length of a shortest cycle in $G$.} {The adjacency matrix $A(G)$ of $G$ is a square $|V| \times |V|$ matrix such that its element $a_{ij}$ is one when there is an edge between vertex $i$ and vertex $j$, and zero when there is no such edge.} For other definitions and notations of graph theory not given here, we refer to \cite{bondy2008graph}.

A set $S \subseteq V$ is called a dominating set if every vertex of $G$ is either in $S$ or adjacent to a vertex of $S$.  The minimum cardinality of a dominating set in $G$, denoted by $\gamma(G)$, is called the domination number.  Domination of graphs has been extensively studied in the scientific literature. The variants of domination have abundant applications, including error-correcting codes constructions for digital communication and efficient data routing in wireless networks \cite{dreyer2000applications,chaovalitwongse2007set,mercier2010survey,cecilio2010survey,castorini2008optimal}. Many different kinds of domination arose, such as the connected dominating set \cite{das1997routing}, the edge dominating set \cite{tural2016maximal}, the total domination \cite{henning2009survey} \cite{macaiyao2019total}, the independent domination \cite{goddard2013independent}, Roman domination \cite{stewart1999defend}, etc.

In this paper, we consider the so-called double Roman domination problem (DRDP), which was initially introduced by Beeler et al. \cite{beeler2016double}. This problem is a generalization of Roman domination problem.

A \textit{Roman dominating function} (RDF) on $G=(V,E)$ is a function $f : V \rightarrow \{0, 1, 2\}$ {with the property that any vertex} $v \in V$ with $f (v) = 0$ is adjacent to at least one vertex $u$ with $f (u) = 2.$ The weight of {an RDF} is defined as $w (f) =\sum_{v \in V}f (v)$. The \textit{Roman domination number} $\gamma_{R}(G)$ is the minimum weight among all RDFs on $G$. {An RDF} on $G$ with weight $\gamma_{R}(G)$ is called a $\gamma_R$-function of $G$. The \textit{Roman domination problem} (RDP) is to determine the value of $\gamma_{R}(G)$.  Roman domination was originally defined and discussed by Stewart et al. \cite{stewart1999defend} , ReVelle and Rosing \cite{revelle2000defendens}, and subsequently developed by Cockayne et al. \cite{cockayne2004roman}. The study of Roman domination was motivated by the defense strategies of the Roman Empire during the reign of Emperor Constantine the Great, $274-337$ AD. To defend the Roman Empire, the emperor decreed that for each city in the Roman Empire, at most two legions should be stationed. Moreover, a location without legion must be adjacent to at least one city at which two legions were stationed, so that if the location was attacked, then one of the two legions could be sent to defend it. So far, rich studies have been performed on various aspects of Roman domination in graphs. For more details, we refer to \cite{ahangar2017maximal,beeler2016double,chambers2009extremal,ivanovic2016improved,chellali2016roman,ahangar2017mixed}.

Inspired by the concept of Roman domination, Beeler et al. \cite{beeler2016double}  proposed a stronger version that doubles the protection by ensuring that any attack can be defended by at least two legions. They allowed at most three legions {to be} stationed at each location. This approach provides a level of defense that is both stronger and more flexible, at less than the anticipated additional cost.

\begin{figure}[htbp]
\centering
\includegraphics[scale=0.32]{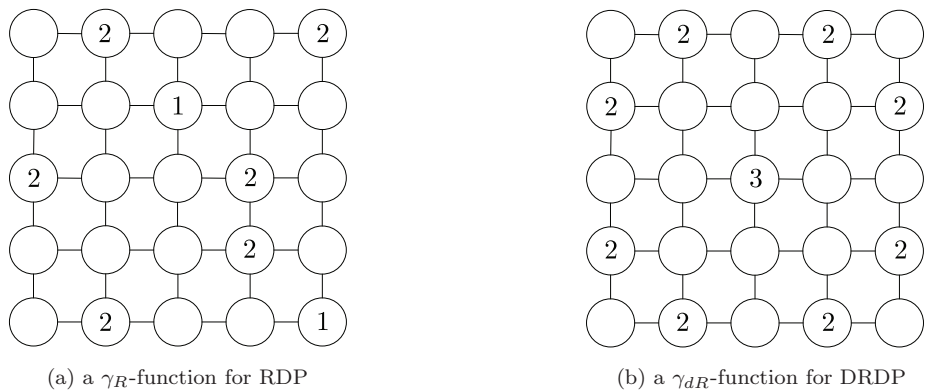}\\
\footnotesize{\vspace*{0.3em} \hspace{0.5em} (a) a $\gamma_R$-function for RDP \hspace{3.9cm} (b) a $\gamma_{dR}$-function for DRDP}\\
\caption{Example of RDP and DRDP on a $5\times5$ grid graph (The dominating function value of the unlabelled vertex is $0$).}
\label{fig}
\end{figure}

Similar to Roman dominating function, a \textit{double Roman dominating function} (DRDF) is a function $f : V \to \{0,1,2,3\}$ having the property that if $f (v) = 0$, then vertex $v$ must have at least two neighbors assigned $2$ under $f$ or {at least} one neighbor $u$ with $f (u) = 3$, and if $f (v) = 1$, then vertex $v$ must have at least one neighbor $u$ with $f (u) \ge 2$. The weight of a DRDF $f$ is the sum $w(f) =\sum_{v\in V} f (v)$, and the \textit{double Roman domination number} of a graph $G$, denoted by $\gamma_{dR}(G)$, is the minimum weight among all possible DRDFs. A DRDF on $G$ with weight $\gamma_{dR}(G)$ is called a $\gamma_{dR}$-function of $G$.
{It is worth noting that Beeler et al. \cite{beeler2016double} showed that for any graph $G$, there exists a $\gamma_{dR}$-function of $G$ such that no vertex needs to be assigned the value $1$.}
The \textit{double Roman domination problem} (DRDP) is to determine the value of $\gamma_{dR}(G)$. An example of the RDP and DRDP is illustrated in Figure \ref{fig}.

The double Roman domination number has been completely determined in paths and cycles  \cite{ahangar2017double}, and has been shown to be linear-time solvable in trees \cite{zhang2018double} and $P_4$-free graphs \cite{yue2018double}. However, for a general graph, it is NP-complete to determine whether a graph has a DRDF of weight at most $k$ \cite{ahangar2017double}. Therefore, it is meaningful to obtain good upper and lower bounds for $\gamma_{dR}(G)$. Beeler et al. \cite{beeler2016double} studied the relationship between $\gamma_{dR}(G)$ and $\gamma(G)$ or $\gamma_{R}(G)$. They showed that $2\gamma(G) \le \gamma_{dR}(G) \le 3\gamma(G)$ and $\gamma_R(G) < \gamma_{dR}(G) \le 2\gamma_{R}(G)$. Moreover, $\gamma_{dR}(T) \ge 2\gamma(T)+1$  for a non-trivial tree $T$. Later, Ahangar et al. \cite{ahangar2017double} gave that a graph with maximum degree $\Delta$ satisfies $\gamma_{dR}(G)\ge \frac{2n}{\Delta} +\frac{\Delta-2}{\Delta}\gamma(G)$. Besides, there are some bounds for $\gamma_{dR}(G)$ in terms of order $n$, minimum degree $\delta$, maximum degree $\Delta$ and diameter $d$. Beeler et al. \cite{beeler2016double} proved that $\gamma_{dR}(G) \le 5n/4$ for a connected graph $G$ with order $n \ge 3$, and characterized the graphs attaining this bound. In \cite{volkmann2018double}, Volkmann obtained the lower bound $\gamma_{dR}(G) \ge \lceil \frac{3n}{\Delta+1} \rceil$ for each nonempty graph $G$. Yue et al. \cite{yue2018double} got the upper bounds $\gamma_{dR}(G) \le 2n -2\Delta +1$ and  $\gamma_{dR}(G) \le \frac{3n(1+\ln{\frac{2(1+\delta)}{3}})}{1+\delta}$. They also showed that almost all graphs have double Roman domination number $\gamma_{dR}\leq n$. Recently,  Anu et al. \cite{anu2018double} found an upper bound with respect to order $n$ and diameter $d$, that is $d+1 \le \gamma_{dR}(G) \le 2n -d$. In \cite{anu2018double}, Anu et al. showed that there is no relation between the double Roman domination number of a graph and its induced subgraph. Yue et al. \cite{yue2018double} obtained that a graph with order $n\ge 3$ satisfies $8 \le \gamma_{dR}(G)+\gamma_{dR}(\overline{G}) \le 2n+3$ and characterized the graphs attaining these bounds. For other results on $\gamma_{dR}(G)$, we refer to \cite{volkmann2018double,jiang2018double,hao2017double,amjadi2018upper,rad2018some}.

However, for a specific graph, {determining the exact value or designing approximation algorithms for the double Roman domination number has not been studied in the literature.} For the RDP, several integer linear programming (ILP) formulations have been proposed by ReVelle and Rosing \cite{revelle2000defendens} and Ivanovi{\'c} \cite{ivanovic2016improved}. Motivated by this, this paper will study the ILP formulations for the DRDP, and will also introduce several extra constraints to strengthen some ILP formulations. Additionally,  an $H(2(\Delta+1))$-approximation algorithm for DRDP is proposed based on one ILP formulation.

The remainder of this paper is organized as follows. In Section \ref{ILP}, we propose two different ILP models for the DRDP. Two improved ILP formulations and two new alternative LP formulations are presented in Section \ref{ALP}. {An approximation algorithm for DRDP is given in Section \ref{Appro-est}.} Finally, we provide computational results in Section \ref{Num-res} and draw our conclusions  in Section \ref{Conclusions}.

\section{ILP models}\label{ILP}
In this section,  we present two different ILP models for the DRDP. These models will be improved in Section \ref{ALP} and experimentally evaluated in Section \ref{Num-res}.

The first ILP model, called \textbf{DRDP-1}, uses three sets of binary variables. Specifically, for each vertex $v \in V$, we define
\begin{equation*}
x_v=\left\{
\begin{array}{ll}
1,&\ f(v)=1\\
0,&\ \text{otherwise}\\
\end{array}
\right.
\hspace{1.5em}
y_v=\left\{
\begin{array}{ll}
1,&\ f(v)=2\\
0,&\ \text{otherwise}\\
\end{array}
\right.
\hspace{1.5em}
z_v=\left\{
\begin{array}{ll}
1,&\ f(v)=3\\
0,&\ \text{otherwise}\\
\end{array}.
\right.
\end{equation*}
Then the DRDP can be formulated as follows:
\begin{subequations}
\begin{align}
\textbf{[DRDP-1]}  \quad  \min \quad & \sum\limits_{v \in V} x_v+2\sum_{v \in V}y_v+3\sum_{v \in V}z_v \label{eq1:1}\\
\text{s.t.}  \quad  & x_v+y_v+z_v+\frac{1}{2}\sum_{u \in N(v)}y_u +\sum_{u \in N(v)}z_u\ge 1,  \quad \forall \; v \in V \label{eq1:2}\\
& \sum_{u \in N(v)}y_u +\sum_{u \in N(v)}z_u\ge x_v, \quad \forall \; v \in V \label{eq1:3}\\
& x_v+y_v+z_v \le 1,  \quad \forall \; v \in V  \label{eq1:4}\\
& x_v,y_v,z_v \in \{0,1\}, \quad \forall \; v \in V \label{eq1:5}
\end{align}
\end{subequations}
The objective function value is given by \eqref{eq1:1}. Constraints \eqref{eq1:2} ensure that if $f (v) = 0$, then vertex $v$ must have at least two neighbors assigned $2$ under $f$ or {at least} one neighbor $u$ with $f (u) = 3$. Constraints \eqref{eq1:3} ensure that if $f (v) = 1$, then vertex $v$ must have at least one neighbor $u$ with $f (u) \ge 2$. By constraints \eqref{eq1:4}, it establishes that every vertex is assigned with at most one value from the set $\{1,2,3\}$.
The decision variables $x_v,\ y_v$ and $z_v$ are set to be binary by constraints \eqref{eq1:5}.

The \textbf{DRDP-1} formulation consists of $3|V|$ binary variables and $3|V|$ constraints. { The following theorem shows that \textbf{DRDP-1} indeed corresponds to the mathematical formulation of the double Roman domination problem.

\begin{thm}\label{equivalence-1}
The optimal objective function value of \textbf{DRDP-1} formulation is equal to the double Roman domination number $\gamma_{dR}$.
\end{thm}

\pf Let $S$ be the set of all the possible double Roman dominating functions and $T$ the set of all feasible solutions of \textbf{DRDP-1} formulation. First we prove that there exists a bijection between $S$ and $T$.

For any double Roman dominating function $f \in S$, we construct a feasible solution $(\textbf{x},\textbf{y},\textbf{z}) \in T \subseteq \{0,1\}^{3\times|V|}$:
 \begin{enumerate}
\item [(a)] if $f(v)=0$, then set $x_v=y_v=z_v=0$;
\item [(b)] if $f(v)=1$, then set $x_v=1, y_v=z_v=0$;
\item [(c)] if $f(v)=2$, then set $y_v=1, x_v=z_v=0$;
\item [(d)] if $f(v)=3$, then set $z_v=1, x_v=y_v=0$.
\end{enumerate}

Note that if $f(v) > 0$, then one of $x_v, y_v, z_v$ equals one, so constraints \eqref{eq1:2} are trivial. Otherwise, $v$ must have at least two neighbors assigned $2$ under $f$ or at least one
neighbor $u$ with $f(u) = 3$, which means there exist at least two vertices $u_1, u_2 \in N(v)$ satisfying $y_{u_1}=y_{u_2}=1$ or at least one vertex $u \in N(v)$ satisfying $z_{u}=1$ from our construction. Hence, constraints \eqref{eq1:2} are also satisfied.

Consider constraints \eqref{eq1:3}. If $f(v) \neq 1$, then $x_v=0$, so constraints \eqref{eq1:3} are trivial. Otherwise, $v$ must
have at least one neighbor $u$ with $f(u) \ge 2$, which means $y_u=1$ or $z_u=1$ from our construction. Therefore, constraints \eqref{eq1:3} hold as well.

Moreover, constraints \eqref{eq1:4}--\eqref{eq1:5} hold naturally from our construction. Consequently, we obtain an injection from $S$ to $T$.

For the other direction, given a feasible solution $(\textbf{x},\textbf{y},\textbf{z}) \in T$, we construct a function $f : V \to \{0,1,2,3\}$ as follows:

\begin{enumerate}
\item [(a$'$)] if $x_v=y_v=z_v=0$, then set $f(v)=0$;
\item [(b$'$)] if $x_v=1, y_v=z_v=0$, then set $f(v)=1$;
\item [(c$'$)] if $y_v=1, x_v=z_v=0$, then set $f(v)=2$;
\item [(d$'$)] if $z_v=1, x_v=y_v=0$, then set $f(v)=3$.
\end{enumerate}

Observe that constraints \eqref{eq1:2} ensure that if $f (v) = 0$, then vertex $v$ must have at least two neighbors assigned $2$ under $f$ or at least one neighbor $u$ with $f (u) = 3$. Constraints \eqref{eq1:3} guarantee that if $f (v) = 1$, then vertex $v$ must have at least one neighbor $u$ with $f (u) \ge 2$. It follows that $f$ is a double Roman dominating function.

We have proved that there exists a bijection between $S$ and $T$. Since the objective function \eqref{eq1:1} corresponds to the weight of a double Roman dominating function, then the optimal objective function value of \textbf{DRDP-1} formulation is equal to the minimum weight among all the double Roman dominating functions, that is $\gamma_{dR}$.

\hfill $\square$

}

The second ILP model, called \textbf{DRDP-2}, follows the same idea as \textbf{DRDP-1}, apart from the definition of some binary variables. Similar as above, for each vertex $v \in V$, we define
\begin{equation*}
p_v=\left\{
\begin{array}{ll}
1,&\ f(v)\ge1\\
0,&\ \text{otherwise}\\
\end{array}
\right.
\hspace{1.5em}
q_v=\left\{
\begin{array}{ll}
1,&\ f(v)\ge2\\
0,&\ \text{otherwise}\\
\end{array}
\right.
\hspace{1.5em}
r_v=\left\{
\begin{array}{ll}
1,&\ f(v)=3\\
0,&\ \text{otherwise}\\
\end{array}.
\right.
\end{equation*}
Then we can obtain another ILP model for DRDP:
\begin{subequations}
\begin{align}
\textbf{[DRDP-2]}  \quad  \min \quad & \sum\limits_{v \in V} p_v+\sum_{v \in V}q_v+\sum_{v \in V}r_v \label{eq2:1}\\
\text{s.t.}  \quad  & p_v+\frac{1}{2}\sum_{u \in N(v)}q_u +\frac{1}{2}\sum_{u \in N(v)}r_u\ge 1,  \quad \forall \; v \in V \label{eq2:2}\\
& q_v +\sum_{u \in N(v)}q_u\ge p_v, \quad \forall \; v \in V \label{eq2:3}\\
& r_v\le q_v\le p_v, \quad \forall \; v \in V  \label{eq2:4}\\
& p_v,q_v,r_v \in \{0,1\}, \quad \forall \; v \in V \label{eq2:5}
\end{align}
\end{subequations}
The objective function value is given by \eqref{eq2:1}. Constraints \eqref{eq2:2} ensure that if $f (v) = 0$, then vertex $v$ must have at least two neighbors assigned $2$ under $f$ or {at least} one neighbor $u$ with $f (u) = 3$. Constraints \eqref{eq2:3} ensure that if $f (v) = 1$, then vertex $v$ must have at least one neighbor $u$ with $f (u) \ge 2$. Constraints \eqref{eq2:4} guarantee that the binary variables are well-defined. Again, the decision variables $p_v,\ q_v$ and $r_v$ are preserved to be binary by constraints \eqref{eq2:5}.

The \textbf{DRDP-2} formulation consists of $3|V|$ binary variables and $4|V|$ constraints. { Similar to the previous arguments, we have the following theorem:

\begin{thm}\label{equivalence-2}
The optimal objective function value of \textbf{DRDP-2} formulation is equal to the double Roman domination number $\gamma_{dR}$.
\end{thm}

}

\section{Alternative ILP models}\label{ALP}
In this section, we first improve the formulations {introduced} in {Section \ref{ILP}}, and then propose two alternative LP formulations.
\subsection{Improved ILP formulations}
To begin with, we state the following result proved by Beeler et al. \cite{beeler2016double}.

\begin{thm}[\cite{beeler2016double}]\label{beeler}
For any graph $G$, there exists a $\gamma_{dR}$-function of $G$ such that no vertex needs to be assigned the value $1$.
\end{thm}

From Theorem \ref{beeler}, we can only consider the $\gamma_{dR}$-functions for which all vertices are assigned the value from the set $\{0,2,3\}$.
Based on this, we substitute $0$ for $x_v$ in \textbf{DRDP-1} to obtain the new ILP model \textbf{DRDP-1$'$}.
\begin{subequations}
\begin{align}
\textbf{[DRDP-1$'$]}  \quad  \min \quad & 2\sum_{v \in V}y_v+3\sum_{v \in V}z_v \label{eq3:1}\\
\text{s.t.}  \quad  & y_v+z_v+\frac{1}{2}\sum_{u \in N(v)}y_u +\sum_{u \in N(v)}z_u\ge 1,  \quad \forall \; v \in V \label{eq3:2}\\
& y_v+z_v \le 1,  \quad \forall \; v \in V  \label{eq3:3}\\
& y_v,z_v \in \{0,1\}, \quad \forall \; v \in V \label{eq3:4}
\end{align}
\end{subequations}

{Using the same idea as the one in Theorem \ref{equivalence-1}, we can obtain the following statement:

\begin{thm}\label{equivalence-3}
The optimal objective function value of \textbf{DRDP-1$'$} formulation is equal to the double Roman domination number $\gamma_{dR}$.
\end{thm}

\pf Let $S$ be the set of all the possible double Roman dominating functions and $T$ the set of all feasible solutions of \textbf{DRDP-1} formulation. Then let $S'$ be the set of all the double Roman dominating functions that no vertex is assigned the value $1$. Clearly, $S' \subseteq S$. Let $T'$ be the set of all feasible solutions of \textbf{DRDP-1$'$} formulation. Note that for every feasible solution $(\textbf{y},\textbf{z}) \in T'$, one can obtain a feasible solution $(\overbrace{0,\ldots,0}^n, \textbf{y},\textbf{z}) \in T$. Since there exists a natural bijection between $T'$ and the set $\{(\overbrace{0,\ldots,0}^n, \textbf{y},\textbf{z}) \ |\ (\textbf{y},\textbf{z}) \in T'\} \subseteq T$,  $T' \subseteq T$ holds. By Theorem \ref{beeler}, $S'$ contains at least one $\gamma_{dR}$-function. It remains to show that there exists a bijection between $S'$ and $T'$ and that the objective function value of \textbf{DRDP-1$'$} corresponds to the weight of a double Roman dominating function in $S'$, whose proof is similar to that in Theorem \ref{equivalence-1} (by excluding the cases $(b)$ and $(b')$).

\hfill $\square$
}

Besides, we can substitute $q_v$ for $p_v$ in \textbf{DRDP-2} by Theorem \ref{beeler}. Then we obtain another improved ILP model \textbf{DRDP-2$'$}.
\begin{subequations}
\begin{align}
\textbf{[DRDP-2$'$]}  \quad  \min \quad & 2\sum_{v \in V}q_v+\sum_{v \in V}r_v \label{eq4:1}\\
\text{s.t.}  \quad  & q_v+\frac{1}{2}\sum_{u \in N(v)}q_u +\frac{1}{2}\sum_{u \in N(v)}r_u\ge 1,  \quad \forall \; v \in V \label{eq4:2}\\
& r_v\le q_v, \quad \forall \; v \in V  \label{eq4:3}\\
& q_v,r_v \in \{0,1\}, \quad \forall \; v \in V \label{eq4:4}
\end{align}
\end{subequations}

{Similar to Theorem \ref{equivalence-3}, we deduce that,

\begin{thm}\label{equivalence-4}
The optimal objective function value of \textbf{DRDP-2$'$} formulation is equal to the double Roman domination number $\gamma_{dR}$.
\end{thm}
}

Next, we explore some extra inequalities to strengthen the formulation \textbf{DRDP-1$'$}. Note that for a $\gamma_{dR}$-function $f$ of graph $G$, the set $\{v\in V \ |\  f(v) \ge 2 \}$ is a dominating set of $G$. Therefore, by the definition of $y_v$ and $z_v$, the sum $\sum_{v \in V}y_v+\sum_{v \in V}z_v$ is equal to the {cardinality} of the set $\{v\in V \ |\  f(v) \ge 2 \}$, and we can {conclude that} the domination number $\gamma(G)$ is bounded by $\sum_{v \in V}y_v+\sum_{v \in V}z_v$. Since $\gamma(G) \ge \lceil n/(1+\Delta) \rceil$ \cite{hedetniemi2013fundamentals}, it is valid to write
\begin{equation*}
\sum_{v \in V}y_v+\sum_{v \in V}z_v \ge \gamma(G) \ge \lceil \frac{n}{1+\Delta} \rceil.
\end{equation*}
Observing that for a $\gamma_{dR}$-function $f$ where all vertices are assigned the values from the set $\{0,2,3\}$, if we replace the function value $3$ by $2$ in $f$, the new function $f'$ is a Roman dominating function (since $f$ is a double Roman dominating function, we know that every vertex $v \in V$ with $f' (v) = 0$ is adjacent to at least one vertex $u$ with $f' (u) = 2$). Thus the weight $2\sum_{v \in V}y_v+2\sum_{v \in V}z_v$ of the function $f'$ is an upper bound for $\gamma_{R}(G)$, that is
\begin{equation*}
2\sum_{v \in V}y_v+2\sum_{v \in V}z_v \ge \gamma_{R}(G).
\end{equation*}
There are some lower bounds for $\gamma_R(G)$ in \cite{mobaraky2008bounds}:
\begin{itemize}
\item [-] $\gamma_{R}(G) \ge \lceil 2g/3 \rceil$, when $G$ contains a cycle;
\item [-] $\gamma_{R}(G) \ge \lceil (d+2)/2 \rceil$, when $G$ is connected;
\item [-] $\gamma_{R}(G) \ge 2\delta$, when $G$ is connected and $g \ge 5$;
\item [-] $\gamma_{R}(G) \ge 4(\delta-1)$, when $G$ is connected, $g \ge 6$ and $\delta \ge 2$;
\item [-] $\gamma_{R}(G) \ge 2\Delta$, when $G$ is connected, $g \ge 7$ and $\delta \ge 2$.
\end{itemize}
Therefore, a valid inequality for \textbf{DRDP-1$'$} is
\begin{equation}\label{L1}
2\sum_{v \in V}y_v+2\sum_{v \in V}z_v \ge L_1, 
\end{equation}
where
\begin{equation*}
L_1 =\max \left\{ 2 \cdot \lceil \frac{n}{1+\Delta} \rceil, \lceil \frac{2g}{3} \rceil \right\}.
\end{equation*}
For convenience, we use $g \ge a$ and $\delta \ge b$ to denote $\{G \ | \ g(G) \ge a\}$ and $\{G \ | \ \delta(G) \ge b\}$, respectively. For a family $\mathcal{A}$ of graphs, we define an indicator function $\textbf{1}_\mathcal{A}$ as follows:
\begin{equation*}
\textbf{1}_\mathcal{A}(G)=\left\{
\begin{array}{ll}
1&\ \text{if} \ G \in \mathcal{A},\\
0&\ \text{if} \ G \notin \mathcal{A}.\\
\end{array}
\right.
\end{equation*}
Further, for a connected graph $G$, we can get a stronger inequality
\begin{equation}\label{L2}
2\sum_{v \in V}y_v+2\sum_{v \in V}z_v \ge L_2 
\end{equation}
where
\begin{equation*}
L_2= \max \left\{ 2 \cdot \lceil \frac{n}{1+\Delta} \rceil, \lceil \frac{d+2}{2} \rceil, \lceil \frac{2g}{3} \rceil, 2\delta \cdot \textbf{1}_{g\ge5},4(\delta-1)\cdot\textbf{1}_{g\ge6 \cap \delta\ge2}, 2\Delta\cdot\textbf{1}_{g\ge7 \cap \delta\ge2} \right\}.
\end{equation*}
Besides, there are some bounds for $\gamma_{dR}(G)$,
which can also be used to derive extra inequalities.
From \cite{anu2018double,volkmann2018double,ahangar2017double}, we know that $\gamma_{dR}(G) \ge \max\{d+1, \lceil 3n/(\Delta+1) \rceil, \lceil 2n/\Delta+(\Delta-2)\gamma(G)/\Delta\rceil\}$. So we can get
\begin{equation}\label{L3}
2\sum_{v \in V}y_v+3\sum_{v \in V}z_v \ge L_3 
\end{equation}
where
\begin{equation*}
L_3=\max \left\{d+1, \lceil \frac{3n}{\Delta+1} \rceil, \lceil \frac{2n}{\Delta}+\frac{\Delta-2}{\Delta}\cdot\frac{n}{1+\Delta} \rceil \right\}.
\end{equation*}
Also, from \cite{ahangar2017double,beeler2016double,anu2018double,yue2018double}, we obtain that
\begin{enumerate}
\item [-] $\gamma_{dR}(G) \le \min \{2n-d, 2n-2\Delta+1, 3n(1+\ln \frac{2(1+\delta)}{3})/(1+\delta) \}$;
\item [-] $\gamma_{dR}(G) \le n$, when $\delta \ge 3$.
\item [-] $\gamma_{dR}(G) \le 5n/4$, when $G$ is a connected graph with order $n \ge 3$.
\end{enumerate}
These upper bounds result in the constraint
\begin{equation}\label{U1}
2\sum_{v \in V}y_v+3\sum_{v \in V}z_v \le U_1 
\end{equation}
where
\begin{equation*}
U_1=\min \left\{2n-d, 2n-2\Delta+1, n\cdot\textbf{1}_{\delta\ge3}+3n\cdot(1-\textbf{1}_{\delta\ge3}), \lfloor \frac{3n(1+\ln \frac{2(1+\delta)}{3})}{1+\delta} \rfloor \right\}.
\end{equation*}
Further, for a connected graph $G$ with order $n \ge 3$,  the above inequality can be strengthened as
\begin{equation}\label{U2}
2\sum_{v \in V}y_v+3\sum_{v \in V}z_v \le U_2 
\end{equation}
where
\begin{equation*}
U_2=\min \left\{2n-d, 2n-2\Delta+1, n\cdot\textbf{1}_{\delta\ge3}+\lfloor \frac{5n}{4} \rfloor \cdot(1-\textbf{1}_{\delta\ge3}), \lfloor \frac{3n(1+\ln \frac{2(1+\delta)}{3})}{1+\delta} \rfloor \right\}.
\end{equation*}

Similar as the above inequalities \eqref{L1}--\eqref{L3} and \eqref{U1}--\eqref{U2}, we can also explore some extra constraints for the formulation \textbf{DRDP-2$'$}. Especially, the sum $2\sum_{v \in V}q_v$ can be viewed as Roman dominating function's weight, which is an upper bound for $\gamma_R(G)$. A series of extra inequalities for formulation \eqref{eq4:1}--\eqref{eq4:4} are listed as follows:
\begin{align}
&2\sum_{v \in V}q_v \ge L_1,\label{L4} \\
&2\sum_{v \in V}q_v+\sum_{v \in V}r_v\ge L_2,\label{L5}\\
&2\sum_{v \in V}q_v+\sum_{v \in V}r_v \le U_1. \label{U3}
\end{align}
In addition, when $G$ is a connected graph with order $n \ge 3$, the above constraints \eqref{L4} and \eqref{U3} can be strengthened as
\begin{equation}\label{L6}
2\sum_{v \in V}q_v \ge L_3, 
\end{equation}
and
\begin{equation}\label{U4}
2\sum_{v \in V}q_v+\sum_{v \in V}r_v \le U_2. 
\end{equation}

We add the constraints \eqref{L1}--\eqref{L3} and \eqref{U1}--\eqref{U2} to \textbf{DRDP-1$'$} and denote the strengthened formulation by \textbf{DRDP-1$'$(+)}. Similarly, \textbf{DRDP-2$'$(+)} stands for the new formulation obtained by adding constraints \eqref{L4}--\eqref{L6} and \eqref{U3}--\eqref{U4} into \textbf{DRDP-2$'$}. In Section \ref{Num-res}, we will explore the computational results for them.

\subsection{Alternatives of the presented ILP formulations}
In this section, {we utilize the similar idea as the one in \cite{ivanovic2016improved} to obtain two alternative LP formulations for DRDP and prove the equivalence between the new formulations and the previous ones. }

The first formulation, denoted by \textbf{DRDP-1$''$}, is obtained from \textbf{DRDP-1$'$} by omitting constraints \eqref{eq3:3}.
\begin{subequations}
\begin{align}
\textbf{[DRDP-1$''$]}  \quad  \min \quad & 2\sum_{v \in V}y_v+3\sum_{v \in V}z_v \label{eq5:1}\\
\text{s.t.}  \quad  & y_v+z_v+\frac{1}{2}\sum_{u \in N(v)}y_u +\sum_{u \in N(v)}z_u\ge 1,  \quad \forall \; v \in V \label{eq5:2}\\
& y_v,z_v \in \{0,1\}, \quad \forall \; v \in V \label{eq5:3}
\end{align}
\end{subequations}

The following theorem proves the equivalence between
\textbf{DRDP-1$''$} and \textbf{DRDP-1$'$}.

\begin{thm}
For any graph $G$, optimal objective function value of  formulation \eqref{eq5:1}--\eqref{eq5:3} is equal to that of  formulation \eqref{eq3:1}--\eqref{eq3:4}.
\end{thm}
\pf Let $V=\{1,\ldots,n\}$. The optimal objective function value of  \textbf{DRDP-1$'$} and \textbf{DRDP-1$''$} is denoted by $OPT_{1'}$ and $OPT_{1''}$, respectively.

First, we prove $OPT_{1''} \ge OPT_{1'}$. Note that an optimal solution to \textbf{DRDP-1$''$} can be viewed as a vector $(\tilde{y}'',\tilde{z}'') \in \{0,1\}^{2n}$, where $\tilde{y}''=(y_1'',\ldots,y_n'')$, $\tilde{z}''=(z_1'',\ldots,z_n'')$. Based on this, we can divide $V$ into two disjoint subsets $V_1$ and $V_2$ such that

\begin{itemize}
\item $V_1=\{\; i \ | \ y_i''=0 \ \text{or} \ z_i''=0\}$,
\item $V_2=\{\; i \ | \ y_i''=1 \ \text{and} \ z_i''=1\}$.
\end{itemize}
Then we define a new vector $(\tilde{y}',\tilde{z}')=(y_1',\ldots,y_n',z_1',\ldots,z_n')$ as follows:

\begin{equation*}
y_i'=\left\{\begin{aligned}
y_i'', \quad & \text{$\forall \; i \in V_1$},\\
0,\ \quad & \text{$\forall \; i \in V_2$},
\end{aligned}\right.~~~
z_i'=z_i'', ~\forall \; i \in V
\end{equation*}
Now, we prove that $(\tilde{y}',\tilde{z}')$ is a feasible solution to \textbf{DRDP-1$'$} formulation.
\begin{itemize}
\item For each $i \in V_1$, by the definition of variables $y_i'$ and $z_i'$ we know that $(y_i',z_i')=(y_i'',z_i'')$, so constraints
\eqref{eq3:2} and \eqref{eq3:4} are satisfied. Further, from the definition of $V_1$, it follows that
\begin{equation*}
y_i'+z_i'=y_i''+z_i''\le 1, \quad \forall \; i \in V_1
\end{equation*}
that is, constraints \eqref{eq3:3} are satisfied.

\item For each $i \in V_2$, by the definition of variables $y_i'$ and $z_i'$, we obtain $(y_i',z_i')=(0,1)$. Then we can verify that constraints \eqref{eq3:2}-\eqref{eq3:4} are satisfied again:
\begin{gather*}
y_i'+z_i'+\frac{1}{2}\sum_{j \in N(i)}y_j' +\sum_{j \in N(i)}z_j'= 0+1+\frac{1}{2}\sum_{j \in N(i)}y_j' +\sum_{j \in N(i)}z_j' \ge 1, \quad \forall \; i \in V_2 \\
y_i'+z_i'=0+1\le 1, \quad \forall \; i \in V_2 \\
y_i'=0\in \{0,1\},\quad z_i'=1 \in \{0,1\}, \quad \forall \; i \in V_2
\end{gather*}
\end{itemize}
Therefore, an optimal solution $(\tilde{y}'',\tilde{z}'')$ to  \textbf{DRDP-1$''$}  corresponds to a feasible solution $(\tilde{y}',\tilde{z}')$ to \textbf{DRDP-1$'$}. It follows that
\begin{align*}
OPT_{1''}&=2\sum_{i \in V}y_i''+3\sum_{i \in V}z_i''=2\sum_{i \in V_1}y_i''+3\sum_{i \in V_1}z_i''+2\sum_{i \in V_2}y_i''+3\sum_{i \in V_2}z_i'' \\
&=2\sum_{i \in V_1}y_i''+3\sum_{i \in V_1}z_i''+5|V_2|,\\
OPT_{1'}& \le 2\sum_{i \in V}y_i'+3\sum_{i \in V}z_i'=2\sum_{i \in V_1}y_i'+3\sum_{i \in V_1}z_i'+2\sum_{i \in V_2}y_i'+3\sum_{i \in V_2}z_i' \\
&=2\sum_{i \in V_1}y_i''+3\sum_{i \in V_1}z_i''+3|V_2|.
\end{align*}
Hence, we get $OPT_{1''} \ge OPT_{1'}$.

It remains to show that $OPT_{1''} \le OPT_{1'}$. Let $S'$, $S''$ be the set of feasible solutions to  \textbf{DRDP-1$'$} and  \textbf{DRDP-1$''$}, respectively.
Note that each feasible solution to \textbf{DRDP-1$'$} is also a feasible solution to \textbf{DRDP-1$''$}. Hence $S' \subseteq S''$, which implies that
\begin{equation*}
OPT_{1''}=\min_{S''} \left(2\sum_{i \in V}y_i+3\sum_{i \in V}z_i\right) \le \min_{S'} \left(2\sum_{i \in V}y_i+3\sum_{i \in V}z_i\right) = OPT_{1'}.
\end{equation*}
This completes our proof.

\hfill $\square$

Next we prove that the binary variables $r_v$ in \textbf{DRDP-2$'$} formulation can be relaxed to non-negative real numbers, and thus we have the following MIP formulation \textbf{DRDP-2$''$}.
\begin{subequations}
\begin{align}
\textbf{[DRDP-2$''$]}  \quad  \min \quad & 2\sum_{v \in V}q_v+\sum_{v \in V}r_v \label{eq6:1}\\
\text{s.t.}  \quad  & q_v+\frac{1}{2}\sum_{u \in N(v)}q_u +\frac{1}{2}\sum_{u \in N(v)}r_u\ge 1,  \quad \forall \; v \in V \label{eq6:2}\\
& r_v\le q_v, \quad \forall \; v \in V  \label{eq6:3}\\
& q_v \in \{0,1\}, \ r_v \in [0, +\infty), \quad \forall \; v \in V \label{eq6:4}
\end{align}
\end{subequations}

\begin{thm}
For any graph $G$, optimal objective function value of formulation \eqref{eq6:1}--\eqref{eq6:4} is equal to that of formulation \eqref{eq4:1}--\eqref{eq4:4}.
\end{thm}
\pf Let $V=\{1,\ldots,n\}$. We denote the optimal objective function value of \textbf{DRDP-2$'$}  and \textbf{DRDP-2$''$} by $OPT_{2'}$ and $OPT_{2''}$, respectively.

Since the optimal objective function value of every relaxed minimization problem is less than or equal to that of the associated original problem, hence we get $OPT_{2''} \le OPT_{2'}$.

Therefore, we only need to prove that $OPT_{2''} \ge OPT_{2'}$. Note that an optimal solution to \textbf{DRDP-2$''$} can be viewed as a vector $(\tilde{q}'',\tilde{r}'') \in \{0,1\}^{2n}$, where $\tilde{q}''=(q_1'',\ldots,q_n'')$, $\tilde{r}''=(r_1'',\ldots,r_n'')$. And we define a new vector $(\tilde{q}',\tilde{r}')=(q_1',\ldots,q_n',r_1',\ldots,r_n')$ as follows:
\begin{equation*}
q_i'=q_i'', \quad \forall i \in V
\end{equation*}
and
\begin{equation*}
r_i'=\left\{\begin{aligned}
0, \quad & \text{ $r_i'' \in [0,1)$}\\
1, \quad & \text{ $r_i'' \in [1,+\infty)$}
\end{aligned}\right.
\end{equation*}

Now we prove that $(\tilde{q}',\tilde{r}')$ is a feasible solution to \textbf{DRDP-2$'$} formulation. The variables $q_i'$ and $r_i'$ have binary values by definition, thus constraints \eqref{eq4:4} hold.
Note that
\begin{enumerate}
\item [i)] when $r_i'' \in [0,1)$, $r_i'=0 \le q_i';$
\item [ii)] when $r_i'' \in [1, +\infty)$, $r_i'=1 \le r_i'' \le q_i'' = q_i'$, where the second inequality follows from constraints \eqref{eq6:3}.
\end{enumerate}
Therefore, constraints \eqref{eq4:3} hold as well. Next we verify the validity of constraints \eqref{eq4:2}. Recall that for each $i \in V$, $q_i''+\frac{1}{2}\sum_{j \in N(i)}q_j'' +\frac{1}{2}\sum_{j \in N(i)}r_j'' \ge 1$ (constraints \ref{eq6:2}) and $r_i'' \le q_i''$ (constraints \ref{eq6:3}). We consider the following three cases:
\begin{enumerate}
\item [i)] $q_i''=1$ or there exist at least two different vertices $j_1, j_2 \in N(i)$ with $q_{j_1}''=q_{j_2}''=1$. Since $q_j'=q_j''$ for each $j \in V$, constraints \eqref{eq4:2} hold naturally for $q_i'$ and $r_i'$.
\item [ii)] $q_i''=0$ and $q_j''=0$ for each $j \in N(i)$. By constraints \eqref{eq6:3}, we have $r_j''=0$ for each $j \in N(i)$, {and hence} $q_i''+\frac{1}{2}\sum_{j \in N(i)}q_j'' +\frac{1}{2}\sum_{j \in N(i)}r_j'' =0$, {which contradicts the constraints \eqref{eq6:2}}. Thus this case is impossible.
\item [iii)] $q_i''=0$ and there exists a unique vertex $j_0\in N(i)$ such that $q_{j_0}''=1$. By constraints \eqref{eq6:3}, we get $\sum_{j \in N(i)\backslash j_0}r_j''=0$. From constraints \eqref{eq6:2}, it follows that $r_{j_0}''\ge 1$. So $q_{j_0}'=r_{j_0}'=1$, then constraints \eqref{eq4:2} are satisfied for $q_i'$ and $r_i'$.
\end{enumerate}

 To sum up, constraints \eqref{eq4:2}--\eqref{eq4:4} are satisfied. Therefore, the vector $(\tilde{q}',\tilde{r}')$ is a feasible solution to \textbf{DRDP-2$'$} formulation, which implies
\begin{equation*}
OPT_{2'} \le 2\sum_{i \in V}q_i'+\sum_{i \in V}r_i' \le  2\sum_{i \in V}q_i''+\sum_{i \in V}r_i''=OPT_{2''}.
\end{equation*}
The second inequality follows from the definition of $(\tilde{q}',\tilde{r}')$.
Thus we complete our proof.

\hfill $\square$\\

\noindent \textbf{Remark.} Actually, from constraints \eqref{eq6:3}, we can  replace constraints \eqref{eq6:4} by
\begin{equation*}
 q_v \in \{0,1\}, \ r_v \in [0, 1], \quad \forall \; v \in V
\end{equation*}


\section{Approximation algorithm}\label{Appro-est}
In this section, we consider the approximation algorithm for the double Roman domination problem. We design a greedy algorithm based on the ILP formulation
\textbf{DRDP-1$''$}. Firstly, consider the following integer programming of general covering problem:
\begin{equation} \label{P}
\begin{aligned}
\min \quad &  \textbf{c}\cdot\textbf{x} \\
\text{s.t.}  \quad  & A\textbf{x} \ge \textbf{b} \\
& \textbf{x} \in \{0,1\}^{2n}
\end{aligned} 
\end{equation}

Let $\textbf{x}=(y_1,\ldots,y_n,z_1,\ldots,z_n)^T$, $\textbf{c}=(\overbrace{2,\ldots,2}^n, \overbrace{3,\ldots,3}^n)$, $\textbf{b}=2\cdot \textbf{1}_n$, $A=\left[
	\begin{array}{cc}
	2I_n+A(G) & 2(I_n+A(G))
	\end{array}
	\right]_{n \times 2n}$,  where $I_n$ is a unit matrix of size $n$. Then we transform the \textbf{DRDP-1$''$} formulation into a constrained general covering problem \eqref{P}. Using the same method described in \cite{dobson1982worst}, we can design a greedy algorithm as follows:


\begin{algorithm}[!h]
	\caption{Greedy} \label{algorithm}
	\begin{algorithmic}[1]
	 \REQUIRE a graph $G$
	 \ENSURE $W_1,W_2$
		\STATE $\textbf{x} \leftarrow \textbf{0},  \ S \leftarrow \{1, \ldots , 2n\} $
		\WHILE {\textbf{b} $\neq$ \textbf{0}}
		\STATE $k \leftarrow \argmin_{j \in S}\{c_j/\sum_{i=1}^n a_{ij}\}$
		\STATE $x_k \leftarrow 1$
		\STATE $S \leftarrow S -\{k\}$
		\FOR{$i=1$ to $n$}
		\STATE $b_i \leftarrow b_i -a_{ik}$
		\FOR{$j \ \text{in}\  S$}
		\STATE $a_{ij} \leftarrow \min\{a_{ij}, b_i\}$
		\ENDFOR
		\ENDFOR
		\ENDWHILE
		\STATE $W_1 \leftarrow \textbf{c}\cdot\textbf{x}$
		\FOR{$i=1$ to $n$}
		\IF{$x_i=x_{i+n}=1$}
		\STATE $x_i \leftarrow 0$
		\ENDIF
		\ENDFOR
		\STATE $W_2 \leftarrow \textbf{c}\cdot\textbf{x}$
		\RETURN $W_1,W_2$
	\end{algorithmic}
\end{algorithm}


\begin{thm}[\cite{dobson1982worst}]\label{dobson}
If OPT is the optimal value of the above covering problem \eqref{P}, $W_1$ is the value given by the greedy algorithm, then
\begin{equation*}
\frac{W_1}{OPT} \le H\left( \max_{1 \le j \le 2n}\sum_{i=1}^n a_{ij} \right),
\end{equation*}
where $H(d)=\sum_{i=1}^d 1/i$ is the first $d$ terms of the harmonic series.
\end{thm}

By using Theorem \ref{dobson}, we can obtain that the greedy algorithm is an $H(2(\Delta+1))$-approximation algorithm for double Roman domination problem.

\begin{cor}\label{factor}
There exists an $H(2(\Delta+1))$-approximation algorithm for double Roman domination problem.
\end{cor}
\pf Let OPT be the optimal objective function value of the formulation \textbf{DRDP-1$''$}. By applying algorithm \ref{algorithm}, we have a feasible solution \textbf{x} satisfying the formulation \textbf{DRDP-1$''$}. Since $W_2 \le W_1$ (step $13$--$19$), with theorem \ref{dobson}, we have
\begin{equation*}
\frac{W_2}{OPT} \le \frac{W_1}{OPT} \le H\left( \max_{1 \le j \le 2n}\sum_{i=1}^n a_{ij} \right) = H(2(\Delta+1)),
\end{equation*}
Thus we complete our proof.
\hfill $\square$\\

Note that $H(d) \le 1+ \ln d$, hence the approximation factor is bound by $O(\ln (2(\Delta+1)))$. Next we will show this approximation factor is best possible in some sense.
Recall that $2\gamma(G)\leq \gamma_{dR}(G)\leq 3\gamma(G)$. We get
the following result:
\begin{lemma}\label{DandDRDP}
\textnormal{(1)} If minimum dominating set problem can be approximated up to a factor of
$\alpha$ in some graph class, then DRDP can be approximated up to a factor of $\frac{3}{2}\alpha$ in the same graph class.

\textnormal{(2)} If DRDP can be approximated up to a factor of $\alpha$ in some graph class, then minimum dominating set problem can be approximated
up to a factor of $\frac{3}{2}\alpha$ in the same graph class.
\end{lemma}
\pf
(1) Let $S$ be the solution returned by some $\alpha$-approximation algorithm for minimum dominating set problem on input $G=(V,E)$. Then $|S|\leq \alpha \gamma(G)$. Let $f$ be the function such that $f(u)=3$ for each $u\in S$ and $f(u)=0$ for each $u\in V\setminus S$. It is easy to see that $f$ is a double Roman dominating function with weight $3|S|$. Thus $3|S|\leq 3\alpha \gamma(G)\leq \frac{3}{2}\alpha \gamma_{dR}(G)$.

(2) Let $f$ be the double Roman dominating function with $(V_0,V_1,V_2,V_3)$ returned by some $\alpha$-approximation algorithm for DRDP on input $G=(V,E)$, where $V_i=\{u\in V\ |\ f(u)=i\}$. Then $|V_1|+2|V_2|+3|V_3|\leq \alpha \gamma_{dR}(G)$.
Obviously, $V_2\cup V_3$ forms a dominating set of $G$. So we get $|V_2|+|V_3|\leq \frac{1}{2}(|V_1|+2|V_2|+3|V_3|)\leq \frac{1}{2} \alpha \gamma_{dR}(G)\leq \frac{3}{2} \alpha \gamma(G)$.

\hfill $\square$
\begin{lemma}[\cite{Chleb2008}]\label{Chleb2008}
Minimum dominating set problem cannot be approximated to within a factor of $(1-\varepsilon) \ln(n)$ in polynomial time for any constant $\varepsilon>0$ unless $NP \subseteq DTIME \left( n^{O( \log ( \log ( n ) ) ) }\right)$. The same results hold also in bipartite and split graphs (hence in chordal graphs, and in complements of chordal graphs as well).
\end{lemma}
Combining the above two lemmas, we can conclude that
\begin{thm}
DRDP cannot be approximated to within a factor of $\frac{2}{3}(1-\varepsilon) \ln(n)$ in polynomial time for any constant $\varepsilon>0$ unless $NP \subseteq DTIME \left( n^{O( \log ( \log ( n ) ) ) }\right)$. The same results hold also in bipartite and split graphs (hence in chordal graphs, and in complements of chordal graphs as well).\hfill $\square$\\
\end{thm}

\noindent \textbf{Remark.} For Roman domination problem, if we use the new ILP formulation presented in \cite{ivanovic2016improved} and set $\textbf{c}=(\overbrace{1,\ldots,1}^n, \overbrace{2,\ldots,2}^n)$,
$A=\left[
	\begin{array}{c c}
	I_n & I_n+A(G)
	\end{array}
\right]_{n \times 2n}$,
	$\textbf{b}=\textbf{1}_n$ in previous constrained general covering problem \eqref{P}, then we can also apply the greedy algorithm \ref{algorithm} to Roman domination problem to obtain the approximate result. Hence we get
\begin{cor}\label{factor}
There exists an $H(\Delta+1)$-approximation algorithm for Roman domination problem.
\end{cor}


\section{Computational results}\label{Num-res}
In this section, computational results which show effectiveness of the proposed DRDP formulations are summarized. All DRDP formulations were solved using CPLEX 12.8 optimization solvers. {We use a Python package --- ``NetworkX''\cite{Hagberg2008exploring} to generate all graphs of different types and sizes and the greedy algorithm was carried out by Python 3.7.0.} All computational experiments were performed on Intel(R) Core(TM) i5-4200U CPU $@$ 1.60GHz 2.30GHz with 8GB RAM under Windows 10 operating system.

For all computational experiments, we used grid graphs, random graphs generated by Erd\H{o}s-R\'enyi $G(n,p)$ model and random trees which were chosen uniformly at random from the set of all trees on $n$ nodes. All the formulations posed above  were applied exactly once to each problem instance. The computation time limit for the applications of CPLEX to each model was $3600$ seconds for each graph. The results are presented in numerical form in three tables: Table \ref{grid} contains the results for all grid graphs. The results for all random graphs and random trees are presented in Table \ref{random} and Table \ref{tree}, respectively. The first {fourth} columns of Table \ref{grid} indicate the name of the instances, the number of vertices ($|V|$), the number of edges ($|E|$) {and the best solution value.} For the Table \ref{random}, the first {fourth} columns are the number of vertices ($|V|$), the number of edges ($|E|$), the edge probability ($p$) to generate random graphs and {the best solution value.} The first {three} columns of Table \ref{tree} indicate the name of the instances, the number of vertices ($|V|$) and {the optimal solution value. Each random graph and random tree} with the same number of vertices has three different instances.  These entries of three tables are followed by the results obtained using all the above formulations and {approximation algorithm.} For each formulation, the computation time required to prove optimality is given in seconds together with the number of nodes visited in the search tree or the number of iterations. {For the greedy algorithm, results and time were presented in the table.} More specifically, {all formulations and the greedy algorithm} are compared with each other and {the best method of DRDP formulations (the least computation time) is marked in bold} in each row. We use the sign ``--'' to show the instances where the optimal solution was not reached within 3600 seconds. {For that in Table \ref{random} and Table \ref{grid}, the best results are shown for the corresponding case when 3600-second limit reaches.}

Table \ref{random} presents the computational results on random graphs. Considering the instances in Table \ref{random}, it can be seen that all six DRDP formulations perform better in dense graphs than in sparse graphs. And we can conclude that the computation time of the six models is closely related to the number of edges in a graph. Besides, comparing the performance of the six models, \textbf{DRDP-1$'$} and \textbf{DRDP-1$''$} have better performance than others in general. Although we have added some extra constraints to \textbf{DRDP-1$'$} and \textbf{DRDP-2$'$}, these new formulations \textbf{DRDP-1$'$}(+) and \textbf{DRDP-2$'$}(+) do not perform very well as we expect.

From Table \ref{grid}, we observe that, on the whole, \textbf{DRDP-1$'$}, \textbf{DRDP-1$'$}(+) and \textbf{DRDP-1$''$} perform better than the other formulations in grid graphs. In particular, for all instances with less than $200$ vertices, \textbf{DRDP-1$''$} outperforms \textbf{DRDP-1$'$}. This was to be expected, as the only difference in the two models is the elimination of a set of constraints. And \textbf{DRDP-1$'$}(+) perform better than \textbf{DRDP-1$'$} in some degree. It can also be observed that \textbf{DRDP-2$'$} and \textbf{DRDP-2$''$} are competitive with \textbf{DRDP-1$'$} and \textbf{DRDP-1$''$} in the context of grid graphs with less than $200$ vertices.

When considering all instances together in Table \ref{tree}, \textbf{DRDP-1$''$} is the best-performing model on random trees, which is followed by \textbf{DRDP-2$'$}. Note that the DRDP problem of tree is much easier than gird or random graph, since the six models can even handle tens of thousands of nodes within the time limit of 3600 seconds. Perhaps it is relevant to the complexity of the problem (One can construct a linear-time algorithm to compute the value of $\gamma_{dR}(T)$ for any tree $T$ \cite{zhang2018double}).

{Furthermore, we can notice that the greedy algorithm was very competitive on grid and random graphs. However, this method performed much worse than any formulation in the case of random trees, especially when the tree's order is very large. The main reason is probably that for random trees, this approximation algorithm does not sufficiently utilize the sparsity of matrix $A$. }

\section{Conclusions}\label{Conclusions}
In this paper, we investigated the double Roman domination problem.
We first proposed two ILP formulations \textbf{DRDP-1} and \textbf{DRDP-2}, {and proved these models indeed correspond to the mathematical formulation of the double Roman domination problem}. By Theorem \ref{beeler}, we can restrict our attention to the functions with no vertex assigned value 1, which leads to two improved models \textbf{DRDP-1$'$} and \textbf{DRDP-2$'$}. Next, we explore some extra inequalities to strengthen the formulations \textbf{DRDP-1$'$} and \textbf{DRDP-2$'$}. Furthermore, it was shown that constraints \eqref{eq3:3} could be excluded from \textbf{DRDP-1$'$} and variables $z_v$ in \textbf{DRDP-2$'$} can be relaxed to non-negative real numbers. Thus we obtained two alternative formulations \textbf{DRDP-1$''$} and \textbf{DRDP-2$''$}.  Finally, we use one ILP formulation \textbf{DRDP-1$''$} to give an $H(2(\Delta+1))$-approximation algorithm and provide an inapproximability result for this problem.

Our computational experiments were carried out on grid graphs and random graphs with up to $370$ vertices, and random trees up to $30000$ vertices. The results showed that \textbf{DRDP-1$'$} and  \textbf{DRDP-1$''$} outperformed the other formulations in most situations. Especially, \textbf{DRDP-1$''$} is the best-performing model on random trees and the formulations with extra constraints did not perform very well like we expect. {We also notice that the greedy algorithm was very competitive on grid and random graphs.} We made several improvements (decrease some variables or add some extra constraints) in some formulations, but significant improvements in computational effort for strengthening those formulations were still expected. Given that optimization solver was not efficient in solving large sparse graphs, designing new exact methods or other approximation algorithms for DRDP is an interesting direction for future work. \\

\noindent {\bf Acknowledgments.} Q. Cai is partially supported by National Natural Science Foundation of China (No. 11701297) and Open Project Foundation of Intelligent Information Processing Key Laboratory of Shanxi Province (No. CICIP2018005). Y. Shi and S. Yao are partially supported by National Natural Science Foundation of China (No. 11771221 and 11811540390), Natural Science Foundation of Tianjin, China (No. 17JCQNJC00300) and China-Slovenia bilateral project ``Some topics in modern graph theory" (No. 12-6).

\begin{landscape}
\begin{table}[h!]
  \centering
  \setlength{\belowcaptionskip}{8pt}
  \caption{Numerical results for grid graphs} \label{random}%
   \resizebox{1.558\textwidth}{5.5em}{
        \begin{tabular}{ccccrcccrcccrcccrcccrcccrcccrcccrcccrcc}
    \toprule
    \multicolumn{4}{c}{Instance} &   & \multicolumn{3}{c}{DRDP-1} &   & \multicolumn{3}{c}{DRDP-1$'$} &   & \multicolumn{3}{c}{DRDP-1$'$(+)} &   & \multicolumn{3}{c}{DRDP-1$''$} &   & \multicolumn{3}{c}{DRDP-2} &   & \multicolumn{3}{c}{DRDP-2$'$} &   & \multicolumn{3}{c}{DRDP-2$'$(+)} &   & \multicolumn{3}{c}{DRDP-2$''$} &   & \multicolumn{2}{c}{Greedy} \\
\cmidrule{1-4}\cmidrule{6-8}\cmidrule{10-12}\cmidrule{14-16}\cmidrule{18-20}\cmidrule{22-24}\cmidrule{26-28}\cmidrule{30-32}\cmidrule{34-36}\cmidrule{38-39}    name & $|V|$ & $|E|$ & \textbf{best} &   & \textbf{result} & \textbf{time} & \textbf{nodes} &   & \textbf{result} & \textbf{time} & \textbf{nodes} &   & \textbf{result} & \textbf{time} & \textbf{nodes} &   & \textbf{result} & \textbf{time} & \textbf{nodes} &   & \textbf{result} & \textbf{time} & \textbf{nodes} &   & \textbf{result} & \textbf{time} & \textbf{nodes} &   & \textbf{result} & \textbf{time} & \textbf{nodes} &   & \textbf{result} & \textbf{time} & \textbf{nodes} &   & \textbf{result} & \textbf{time} \\
    \midrule
    Grid05$\times$10 & 50 & 85 & 38 &   & 38 & 0.63 & 79 &   & \textbf{38} & \textbf{0.20} & \textbf{0} &   & 38 & 0.31  & 0 &   & 38 & 0.41 & 13 &   & 38 & 0.81 & 81 &   & 38 & 0.31 & 0 &   & 38 & 0.41  & 0 &   & 38 & 0.63  & 163 &   & 54 & 0.39  \\
    Grid05$\times$15 & 75 & 130 & 56 &   & 56 & 4.06 & 2528 &   & 56 & 1.20  & 630 &   & 56 & 1.73  & 550 &   & \textbf{56} & \textbf{0.72} & \textbf{201} &   & 56 & 1.7 & 2236 &   & 56 & 0.81 & 169 &   & 56 & 1.80  & 232 &   & 56 & 1.89  & 2821 &   & 80 & 1.07  \\
    Grid05$\times$20 & 100 & 175 & 74 &   & 74 & 29.44 & 15505 &   & 74 & 4.86  & 3348 &   & 74 & 3.41  & 2105 &   & 74 & 3.31 & 2576 &   & 74 & 38.48 & 56261 &   & \textbf{74} & \textbf{2.67} & \textbf{2773} &   & 74 & 3.02  & 2607 &   & 74 & 6.14  & 7440 &   & 105 & 2.19  \\
    Grid05$\times$25 & 125 & 220 & 92 &   & 92 & 200.39 & 87272 &   & 92 & 15.69  & 9631 &   & 92 & 14.50  & 5501 &   & \textbf{92} & \textbf{10.22} & \textbf{5210} &   & 92 & 173.11 & 89529 &   & 92 & 17.88 & 15051 &   & 92 & 32.27  & 29349 &   & 92 & 30.20  & 19124 &   & 134 & 4.21  \\
    Grid05$\times$30 & 150 & 265 & 110 &   & 110 & -- & -- &   & 110 & 235.53 & 73607 &   & 110 & 212.36  & 108141 &   & \textbf{110} & \textbf{135.11} & \textbf{84417} &   & 111 & -- & -- &   & 110 & 259.31 & 124201 &   & 110 & 220.97  & 115958 &   & 110 & 293.45  & 147695 &   & 160 & 7.35  \\
    Grid05$\times$35 & 175 & 310 & 128 &   & 129 & -- & -- &   & 128 & 507.59  & 142150 &   & \textbf{128} & \textbf{317.27} & \textbf{112420} &   & 128 & 353.78 & 92031 &   & 129 & -- & -- &   & 128 & 583.75 & 141225 &   & 128 & 537.42  & 178804 &   & 128 & 894.69  & 362089 &   & 185 & 12.33  \\
    Grid10$\times$10 & 100 & 180 & 72 &   & 72 & 22.94 & 18817 &   & 72 & 2.45  & 1418 &   & 72 & 2.03  & 796 &   & \textbf{72} & \textbf{1.25} & \textbf{461} &   & 72 & 16.19 & 14135 &   & 72 & 2.28 & 2335 &   & 72 & 2.56  & 2965 &   & 72 & 6.00  & 5554 &   & 90 & 1.78  \\
    Grid10$\times$15 & 150 & 275 & 106 &   & 106 & 738.36 & 140939 &   & 106 & 89.03  & 31426 &   & 106 & 76.48  & 32497 &   & \textbf{106} & \textbf{21.44} & \textbf{9335} &   & 106 & 493.26 & 160496 &   & 106 & 65.92 & 26816 &   & 106 & 180.01 & 144437 &   & 106 & 97.27  & 40352 &   & 138 & 6.10  \\
    Grid10$\times$20 & 200 & 370 & 140 &   & 141 & -- & -- &   & \textbf{140} & \textbf{1180.09} & \textbf{215452} &   & 140 & 2014.41  & 325629 &   & 140 & 3362.61  & 619175 &   & 140 & -- & -- &   & 141 & -- & -- &   & 140 & 2744.03  & 513429 &   & 139.99 & -- & -- &   & 182 & 13.98  \\
    Grid10$\times$25 & 250 & 465 & 174 &   & 175 & -- & -- &   & 174 & -- & -- &   & 174 & -- & -- &   & 175 & -- & -- &   & 174 & -- & -- &   & 174 & -- & -- &   & 175 & -- & -- &   & 175 & -- & -- &   & 220 & 25.90  \\
    Grid10$\times$30 & 300 & 560 & 208 &   & 209 & -- & -- &   & 211 & -- & -- &   & 209 & -- & -- &   & 208 & -- & -- &   & 208 & -- & -- &   & 209 & -- & -- &   & 210 & -- & -- &   & 210 & -- & -- &   & 268 & 44.43  \\
    Grid15$\times$15 & 225 & 420 & 155 &   & 155 & -- & -- &   & 155 & 1373.06  & 215668 &   & \textbf{155} & \textbf{1205.52} & \textbf{168916} &   & 155 & 1355.20  & 192572 &   & 156 & -- & -- &   & 155 & 1709.80  & 225916 &   & 155 & 2843.16  & 341775 &   & 155 & 1905.75  & 392388 &   & 194 & 19.21  \\
    Grid15$\times$20 & 300 & 565 & 204 &   & 205 & -- & -- &   & 204 & -- & -- &   & 204 & -- & -- &   & 204 & -- & -- &   & 205 & -- & -- &   & 205 & -- & -- &   & 205 & -- & -- &   & 204.99 & -- & -- &   & 258 & 42.19  \\
    Grid15$\times$25 & 375 & 710 & 254 &   & 254 & -- & -- &   & 254 & -- & -- &   & 259 & -- & -- &   & 254 & -- & -- &   & 253 & -- & -- &   & 254 & -- & -- &   & 254 & -- & -- &   & 254 & -- & -- &   & 323 & 79.26  \\
    \bottomrule
    \end{tabular}%
}
\end{table}%

\begin{table}[h!]
  \centering
  \setlength{\belowcaptionskip}{8pt}
  \caption{Numerical results for random graphs}  \label{grid}
  \resizebox{1.558\textwidth}{6em}{
    \begin{tabular}{ccccrcccrcccrcccrcccrcccrcccrcccrcccrcc}
    \toprule
    \multicolumn{4}{c}{Instance} &   & \multicolumn{3}{c}{DRDP-1} &   & \multicolumn{3}{c}{DRDP-1$'$} &   & \multicolumn{3}{c}{DRDP-1$'$(+)} &   & \multicolumn{3}{c}{DRDP-1$''$} &   & \multicolumn{3}{c}{DRDP-2} &   & \multicolumn{3}{c}{DRDP-2$'$} &   & \multicolumn{3}{c}{DRDP-2$'$(+)} &   & \multicolumn{3}{c}{DRDP-2$''$} &   & \multicolumn{2}{c}{Greedy} \\
\cmidrule{1-4}\cmidrule{6-8}\cmidrule{10-12}\cmidrule{14-16}\cmidrule{18-20}\cmidrule{22-24}\cmidrule{26-28}\cmidrule{30-32}\cmidrule{34-36}\cmidrule{38-39}    $|V|$ & $|E|$ & $p$ & \textbf{best} &   & \textbf{result} & \textbf{time} & \textbf{nodes} &   & \textbf{result} & \textbf{time} & \textbf{nodes} &   & \textbf{result} & \textbf{time} & \textbf{nodes} &   & \textbf{result} & \textbf{time} & \textbf{nodes} &   & \textbf{result} & \textbf{time} & \textbf{nodes} &   & \textbf{result} & \textbf{time} & \textbf{nodes} &   & \textbf{result} & \textbf{time} & \textbf{nodes} &   & \textbf{result} & \textbf{time} & \textbf{nodes} &   & \textbf{result} & \textbf{time} \\
    \midrule
    \multirow{3}[2]{*}{100} & 994 & 0.2 & 21 &   & 21 & 3.80  & 5724 &   & \textbf{21} & \textbf{2.22} & \textbf{5162} &   & 21 & 3.30  & 4827 &   & 21 & 2.42  & 5427 &   & 21 & 15.38  & 16799 &   & 21 & 6.08  & 8283 &   & 21 & 3.05  & 4994 &   & 21 & 19.52  & 32246 &   & 23 & 0.40  \\
      & 2484 & 0.5 & 11 &   & 11 & 7.23  & 7519 &   & \textbf{11} & \textbf{4.45} & \textbf{8344} &   & 11 & 5.58  & 8356 &   & 11 & 4.53  & 8676 &   & 11 & 30.38  & 41450 &   & 11 & 6.70  & 8778 &   & 11 & 6.48  & 9196 &   & 11 & 256.86  & 302139 &   & 12 & 0.26  \\
      & 3949 & 0.8 & 6 &   & 6 & 0.55  & 197 &   & \textbf{6} & \textbf{0.53}  & \textbf{0} &   & 6 & 1.16  & 185 &   & 6 & 0.75  & 181 &   & 6 & 0.97  & 201 &   & 6 & 1.09  & 168 &   & 6 & 1.17  & 189 &   & 6 & 4.30  & 8249 &   & 6 & 0.17 \\
    \midrule
    \multirow{3}[2]{*}{150} & 2239 & 0.2 & 24 &   & 24 & 535.69  & 550441 &   & 24 & 1745.03  & 1987548 &   & 24 & 429.39  & 477868 &   & 24 & 3506.45  & 3619108 &   & 26 & -- & -- &   & 24 & 788.81  & 1176644 &   & \textbf{24} & \textbf{402.88} & \textbf{474078} &   & 23.99 & -- & -- &   & 27 & 0.89  \\
      & 5596 & 0.5 & 11 &   & 11 & 73.67  & 54703 &   & 11 & 45.56  & 57470 &   & 11 & 31.25  & 31104 &   & 11 & 39.91  & 56138 &   & 11 & 422.19  & 338528 &   & 11 & 69.92  & 64322 &   & \textbf{11} & \textbf{29.16} & \textbf{18769} &   & 11 & 1523.83  & 2121825 &   & 12 & 0.45  \\
      & 8903 & 0.8 & 6 &   & 6 & 1.97  & 298 &   & 6 & 1.56  & 302 &   & 6 & 2.19  & 291 &   & \textbf{6} & \textbf{1.44}  & \textbf{298} &   & 6 & 2.83  & 221 &   & 6 & 2.39  & 267 &   & 6 & 2.44  & 212 &   & 6 & 16.48  & 12621 &   & 6 & 0.25 \\
    \midrule
    \multirow{3}[2]{*}{200} & 0 & 0.2 & 27 &   & 27 & -- & -- &   & 27 & -- & -- &   & 27 & -- & -- &   & 27 & -- & -- &   & 27 & -- & -- &   & 27 & -- & -- &   & 27 & -- & -- &   & 26.99 & -- & -- &   & 29 & 1.56  \\
      & 9941 & 0.5 & 12 &   & 12 & 219.22  & 100591 &   & 12 & 139.25  & 93823 &   & 12 & 216.39  & 123961 &   & \textbf{12} & \textbf{122.53} & \textbf{100557} &   & 12 & 2159.20  & 986301 &   & 12 & 225.19  & 135567 &   & 12 & 239.94  & 107938 &   & 12 & -- & -- &   & 14 & 0.86  \\
      & 15850 & 0.8 & 6 &   & 6 & 3.95  & 387 &   & \textbf{6} & \textbf{3.16} & \textbf{397} &   & 6 & 5.47  & 427 &   & 6 & 3.84  & 395 &   & 6 & 8.23  & 312 &   & 6 & 5.44  & 381 &   & 6 & 5.16  & 287 &   & 6 & 46.02  & 20269 &   & 9 & 0.58  \\
    \midrule
    \multirow{3}[2]{*}{250} & 6205 & 0.2 & 29 &   & 29 & -- & -- &   & 29 & -- & -- &   & 30 & -- & -- &   & 29 & -- & -- &   & 29 & -- & -- &   & 30 & -- & -- &   & 29 & -- & -- &   & 29 & -- & -- &   & 31 & 2.81  \\
      & 15527 & 0.5 & 12 &   & 12 & 479.44  & 127117 &   & 12 & 294.11  & 122411 &   & 12 & 295.64  & 125937 &   & \textbf{12} & \textbf{253.73}  & \textbf{127877} &   & 12 & -- & -- &   & 12 & 551.56  & 174352 &   & 12 & 500.33  & 143796 &   & 12 & -- & -- &   & 12 & 1.07 \\
      & 24869 & 0.8 & 6 &   & 6 & 31.13  & 1906 &   & 6 & 9.83  & 2729 &   & 6 & 14.86  & 1743 &   & \textbf{6} & \textbf{7.08} & \textbf{995} &   & 6 & 70.42  & 3704 &   & 6 & 12.45  & 493 &   & 6 & 46.98  & 3059 &   & 6 & 142.36  & 37039 &   & 8 & 0.80  \\
    \midrule
    \multirow{3}[2]{*}{300} & 8948 & 0.2 & 30 &   & 30 & -- & -- &   & 30 & -- & -- &   & 30 & -- & -- &   & 30 & -- & -- &   & 30 & -- & -- &   & 30 & -- & -- &   & 30 & -- & -- &   & 30 & -- & -- &   & 30 & 4.51  \\
      & 22393 & 0.5 & 12 &   & 12 & 1499.45  & 244786 &   & \textbf{12} & \textbf{802.33} & \textbf{217587} &   & 12 & 1319.98  & 349369 &   & 14 & -- & -- &   & 14 & -- & -- &   & 12 & 1470.81  & 252751 &   & 14 & -- & -- &   & 12 & -- & -- &   & 15 & 1.97  \\
      & 35902 & 0.8 & 8 &   & 8 & 270.17  & 57352.00  &   & 8 & 113.47  & 37007 &   & 8 & 148.28  & 34773 &   & \textbf{8} & \textbf{99.38} & \textbf{37041} &   & 8 & 1289.20  & 128757 &   & 8 & 180.73  & 34821 &   & 8 & 223.06  & 31940 &   & 8 & -- & -- &   & 10 & 1.46  \\
    \bottomrule
    \end{tabular}%
}
\end{table}%

\begin{table}[h!]
  \centering
  \setlength{\belowcaptionskip}{8pt}
  \caption{Numerical results for random trees}  \label{tree}
  \resizebox{1.6\textwidth}{8em}{
    \begin{tabular}{cccccccccccccccccccccccccccccc}
    \toprule
    \multicolumn{3}{c}{Instance} &   & \multicolumn{2}{c}{DRDP-1} &   & \multicolumn{2}{c}{DRDP-1$'$} &   & \multicolumn{2}{c}{DRDP-1$'$(+)} &   & \multicolumn{2}{c}{DRDP-1$''$} &   & \multicolumn{2}{c}{DRDP-2} &   & \multicolumn{2}{c}{DRDP-2$'$} &   & \multicolumn{2}{c}{DRDP-2$'$(+)} &   & \multicolumn{2}{c}{DRDP-2$''$} &   & \multicolumn{2}{c}{Greedy} \\
\cmidrule{1-3}\cmidrule{5-6}\cmidrule{8-9}\cmidrule{11-12}\cmidrule{14-15}\cmidrule{17-18}\cmidrule{20-21}\cmidrule{23-24}\cmidrule{26-27}\cmidrule{29-30}    name & $|V|$ & \textbf{value} &   & \textbf{time} & \textbf{iterations} &   & \textbf{time} & \textbf{iterations} &   & \textbf{time} & \textbf{iterations} &   & \textbf{time} & \textbf{iterations} &   & \textbf{time} & \textbf{iterations} &   & \textbf{time} & \textbf{iterations} &   & \textbf{time} & \textbf{iterations} &   & \textbf{time} & \textbf{iterations} &   & \textbf{time} & \textbf{result} \\
    \midrule
    Tree-100-1 & \multirow{3}[1]{*}{100} & 102 &   & 0.22  & 371 &   & 0.14  & 185 &   & 0.08  & 198 &   & \textbf{0.03} & \textbf{17} &   & 0.41  & 371 &   & 0.05  & 182 &   & 0.30  & 195 &   & 0.06  & 165 &   & 1.78 & 113 \\
    Tree-100-2 &   & 100 &   & 0.19  & 347 &   & 0.17  & 179 &   & 0.27  & 195 &   & \textbf{0.05} & \textbf{5} &   & 0.55  & 374 &   & 0.09  & 167 &   & 0.25  & 207 &   & 0.08  & 166 &   & 1.92 & 114 \\
    Tree-100-3 &   & 98 &   & 0.16  & 269 &   & 0.06  & 152 &   & 0.11  & 167 &   & \textbf{0.02} & \textbf{3} &   & 0.33  & 280 &   & 0.08  & 136 &   & 0.05  & 187 &   & 0.05  & 132 &   & 1.74 & 110 \\
	\midrule
    Tree-1000-1 & \multirow{3}[1]{*}{1000} & 998 &   & 2.20  & 3258 &   & 0.52  & 1778 &   & 0.75  & 1926 &   & \textbf{0.11} & \textbf{20} &   & 1.55  & 3480 &   & 0.42  & 1654 &   & 0.47  & 1898 &   & 0.61  & 1576 &   & 1681.76 & 1149 \\
   Tree-1000-2 &   & 964 &   & 1.11  & 2993 &   & 0.78  & 1673 &   & 0.86  & 1774 &   & \textbf{0.08} & \textbf{16} &   & 1.16  & 3132 &   & 0.50  & 903 &   & 0.63  & 1706 &   & 0.88  & 1167 &   & 1523.83 & 1092 \\
   Tree-1000-3 &   & 982 &   & 1.30  & 3224 &   & 0.33  & 1719 &   & 0.80  & 1779 &   & \textbf{0.25} & \textbf{44} &   & 1.44  & 3318 &   & 0.53  & 1395 &   & 0.97  & 1763 &   & 0.73  & 1342 &   & 1628.42 & 1136 \\
	\midrule
    Tree-10000-1 & \multirow{3}[1]{*}{10000} & 9838 &   & 15.84  & 32244 &   & 5.17  & 18151 &   & 7.25  & 1913 &   & \textbf{0.66} & \textbf{29} &   & 73.94  & 66361 &   & 3.17  & 14733 &   & 4.45  & 710 &   & 5.11  & 15042 &   & -- & -- \\
    Tree-10000-2 &   & 9833 &   & 13.73  & 32843 &   & 10.25  & 18082 &   & 7.44  & 1982 &   & \textbf{2.14} & \textbf{60} &   & 82.91  & 63080 &   & 3.52  & 15539 &   & 4.59  & 665 &   & 6.73  & 15183 &   & -- & -- \\
    Tree-10000-3 &   & 9868 &   & 26.83  & 32579 &   & 7.78  & 18064 &   & 7.09  & 2169 &   & \textbf{2.39} & \textbf{59} &   & 101.69  & 67826 &   & 5.28  & 15391 &   & 7.13  & 663 &   & 5.64  & 15166 &   & -- & -- \\
	\midrule
    Tree-20000-1 & \multirow{3}[1]{*}{20000} & 19624 &   & 33.16  & 22924 &   & 7.81  & 4063 &   & 11.58  & 3807 &   & \textbf{1.16} & \textbf{123} &   & 97.83  & 70320 &   & 5.75  & 1317 &   & 10.23  & 1325 &   & 7.91  & 29821 &   & -- & -- \\
    Tree-20000-2 &   & 19622 &   & 30.64  & 22741 &   & 10.03  & 4205 &   & 13.83  & 3975 &   & \textbf{1.7} & \textbf{61} &   & 103.05  & 71581 &   & 5.61  & 1398 &   & 10.63  & 1510 &   & 8.03  & 30133 &   & -- & -- \\
    Tree-20000-3 &   & 19655 &   & 34.63  & 22315 &   & 8.08  & 4101 &   & 18.59  & 3924 &   & \textbf{3.19} & \textbf{828} &   & 95.89  & 70016 &   & 5.86  & 1352 &   & 13.05  & 1425 &   & 9.25  & 30120 &   & -- & -- \\
    \midrule
    Tree-30000-1 & \multirow{3}[1]{*}{30000} & 29410 &   & 47.25  & 32731 &   & 12.08  & 5894 &   & 23.11  & 5817 &   & \textbf{2.25} & \textbf{393} &   & 161.22  & 103441 &   & 9.00  & 1996 &   & 19.27  & 2116 &   & 13.27  & 5500 &   & -- & -- \\
    Tree-30000-2 &   & 29309 &   & 44.33  & 31951 &   & 13.47  & 5816 &   & 21.70  & 5575 &   & \textbf{2.09} & \textbf{58} &   & 161.00  & 103398 &   & 9.58  & 2123 &   & 18.97  & 2225 &   & 13.52  & 5302 &   & -- & -- \\
    Tree-30000-3 &   & 29423 &   & 44.56  & 32812 &   & 19.19  & 5961 &   & 21.98  & 6098 &   & \textbf{2.27} & \textbf{458} &   & 149.20  & 102554 &   & 9.77  & 2024 &   & 17.56  & 2093 &   & 13.51  & 5532 &   & -- & -- \\
    \bottomrule
    \end{tabular}%
}
\end{table}%
\end{landscape}


\end{document}